\documentclass[10pt,reqno]{amsart}

\usepackage{a4wide}
\usepackage{amssymb}
\usepackage[all]{xy}
\usepackage[mathscr]{eucal}

\newtheorem{prop}{Proposition}[section]
\newtheorem{lem}[prop]{Lemma}
\newtheorem{thm}[prop]{Theorem}
\newtheorem{cor}[prop]{Corollary}

\theoremstyle{remark}
\newtheorem{rem}[prop]{Remark}

\theoremstyle{definition}
\newtheorem{defn}[prop]{Definition}

\numberwithin{equation}{section}
\numberwithin{prop}{section}

\newcommand{\CC}{\mathbb{C}}
\newcommand{\cst}{\mathrm{C^*}}
\newcommand{\id}{\mathrm{id}}
\newcommand{\comp}{\!\circ\!}
\newcommand{\tens}{\otimes}
\newcommand{\vt}{\!\vartriangle\!}
\newcommand{\Eta}{\boldsymbol{\eta}}

\newcommand{\gC}{\mathfrak{C}}
\newcommand{\gQ}{\mathfrak{Q}}

\newcommand{\st}{\,\vline\:}
\newcommand{\tp}{\xymatrix{*+<.7ex>[o][F-]{\scriptstyle\top}}}
\DeclareMathOperator{\Mor}{Mor}
\DeclareMathOperator{\M}{M}
\DeclareMathOperator{\qs}{\mathscr{QS}}
\DeclareMathOperator{\QMap}{\mathscr{Q}-Map}
\DeclareMathOperator{\Map}{Map}
\DeclareMathOperator{\C}{C}
\DeclareMathOperator{\Cinf}{\C_\infty}

\renewcommand{\Bar}[1]{\overline{#1}}

\begin{document}

\title{Quantum families of maps and quantum semigroups
on finite quantum spaces}

% \date{October 30, 2006}
\date{\today}

\author{Piotr Miko{\l}aj So{\l}tan}
\address{Department of Mathematical Methods in Physics,
Faculty of Physics, Warsaw University}
\email{piotr.soltan@fuw.edu.pl}

\thanks{Research partially supported by KBN grants nos.~2P03A04022
\& 115/E-343/SPB/6.PRUE/DIE50/2005-2008.}

\begin{abstract}
Quantum families of maps between quantum spaces are defined and studied. We
prove that quantum semigroup (and sometimes quantum group) structures arise
naturally on such objects out of more fundamental properties. As particular
cases we study quantum semigroups of maps preserving a fixed state and quantum
commutants of given quantum families of maps.
\end{abstract}

\maketitle

\section{Introduction}

Let $X$ be a set. Then the set $\Map(X)$ of all maps $X\to{X}$ is a
semigroup. Of course, the set of all maps fixing a given point of $X$ is a
subsemigroup of $\Map(X)$. So is the set of all maps leaving invariant a given
measure on $X$ or commuting with a fixed family of maps $X\to{X}$. These
statements border triviality. The situation does not change substantially if we
introduce a topology on $X$ and require that all maps be continuous.

Our aim in this paper is to investigate the non commutative analogs of the
above mentioned phenomena. More precisely we will recall the definition of a
quantum family of maps between quantum spaces (\cite{pseu}) and we will show
that, just as in the commutative case mentioned above, quantum semigroup
structures appear naturally on many such objects.

In \cite{Qsym} S.~Wang investigated quantum automorphism groups of finite
quantum spaces. He searched for universal objects in the category of quantum
transformation groups of a given finite quantum space. He also mentioned
quantum semigroups in \cite[Remark (3), page 208]{Qsym}. We show that the
quantum semigroup structure on these objects is there as a consequence of a
more fundamental structure these objects possess. This makes them also easier
to define.

Wang proved that for non classical finite quantum spaces the quantum
automorphism group does not exist and turned to study the group preserving a
fixed state. We will take a more general approach and define the
quantum family of all maps preserving a given state. Again the quantum
semigroup structure (and in some cases quantum \emph{group} structure ---
cf.~Section \ref{gr}) appears from a more fundamental property of these
objects.

We will give one more example of a similar situation, where a quantum
subsemigroup of a given quantum semigroup is defined without reference to its
semigroup structure, by constructing the quantum commutant of a given
quantum family of maps. Again, the emergence of a quantum seigroup structure
will be a consequence of a more fundamental property of the quantum commutant.

Let us now briefly describe the contents of the paper. Section \ref{qs} is a
short summary of the standard language of non commutative topology. In
particular we shall recall and discuss the definition of a \emph{quantum space}.
 In Section \ref{qspm} we shall define the concept of a quantum family of maps
from one quantum space to another. This notion was introduced already in
\cite{pseu}, where quantum spaces were called ``pseudospaces''. We shall define
what the \emph{quantum space of all maps} from one quantum space to another is
and prove its existence in a special (yet interesting) case. Then we shall
define the crucial notion of composition of quantum families of maps.

The quantum space of all maps from a given finite quantum space to itself
carries a natural structure of a compact quantum semigroup with unit. This is
the content of Section \ref{qss}. We shall study the properties of this
quantum semigroup and its action on the finite quantum space, like ergodicity.

The next two sections are devoted to natural constructions of quantum
subsemigroups of the quantum semigroup defined in Section \ref{qss}. First, in
Section \ref{invS}, we define  and study the quantum family of all maps
preserving a fixed state. This family is naturally endowed with a compact
quantum semigroup structure. The existence of this structure follows from
simple considerations concerning composition of quantum families (as defined in
Section \ref{qspm}). Then, in Section \ref{qc}, we construct the quatum
commutant of a given quantum family of maps. This construction bears many
similarities to the one presented in Section \ref{invS}. It is based on the
notion of \emph{commuting quantum families of maps} which is briefly
investigated at the beginning of the section. The quantum commutant has a
natural structure of a quantum semigroup.

The constructions presented in Sections \ref{invS} and \ref{qc} clarify the
mechanism of obtaining quantum semigroup structure. This has never been
investigated before. In the last section we address the question when
quantum \emph{group} structures appear and when they should be expected and show how S.~Wang's quantum automorphism groups of finite spaces (\cite{Qsym}) fit into the framework developed in Section \ref{invS}.

\section{Quantum spaces}\label{qs}

Let $\gC$ be the category of $\cst$-algebras described and studied in
\cite{pseu,gen}. The
objects of $\gC$ are all $\cst$-algebras and for any two objects $A$ and $B$ of
$\gC$ the set $\Mor(A,B)$ consists of all non degenerate $*$-homomorphisms from
$A$ to $\M(B)$. The category $\gQ$ of quantum spaces is by definition the dual
category of $\gC$. By definition the class of objects of $\gQ$ is the same as
the class of objects of $\gC$. Nevertheless for any $\cst$-algebra $A$ we shall
write $\qs(A)$ for $A$ regarded as an object of $\gQ$.

From the point of non commutative geometry (topology) it is natural to
work with objects of $\gQ$. On the other hand all the tools at our disposal are
from the world of $\cst$-algebras. We shall try to introduce a compromise
between the two conventions by declaring that the phrase ``let $\qs(A)$ be a
quantum space'' be taken to mean ``let $A$ be a $\cst$-algebra''. Moreover we
shall use interchangeably the notation $\Phi\in\Mor(A,B)$ and
$\Phi:\qs(B)\to\qs(A)$.

The category of locally compact topological spaces with continuous maps is a
full subcategory of $\gQ$. A quantum space $\qs(A)$ is a locally compact space
if and only if $A$ is a commutative $\cst$-algebra. In this case
$A=\Cinf\bigl(\qs(A)\bigr)$. In widely accepted terminology quantum spaces
corresponding to commutative $\cst$-algebras are called \emph{classical}
spaces.

Many notions from topology are often generalized to the non commutative
setting. As an example let us mention the fact that a quantum space $\qs(A)$ is
called \emph{compact} if $A$ is unital. If $A$ is finite dimensional then
$\qs(A)$ is said to be \emph{finite quantum space}. A more
controversial idea to call a quantum space $\qs(A)$ a \emph{finite dimensional}
if $A$ is finitely generated was proposed in \cite{pseu}.

An interesting step towards a better understanding of the category $\gQ$ was
taken in \cite{dualth,gen} (see also \cite{k-slw}). The results of these papers
show that any (separable) $\cst$-algebra is of the form $\Cinf(\mathbb{X})$,
where $\mathbb{X}$ is a certain $\mathrm{W}^*$-category and $\Cinf(\,\cdot\,)$
has a whole new meaning (which reduces to the old one for commutative
$\cst$-algebras). This means that, despite technical complications, it
is possible to realize quantum spaces as concrete mathematical objects.

\section{Quantum families of maps}\label{qspm}

In this section we shall introduce the objects of our study. These will be
quantum spaces of maps or quantum families of maps. The latter concept is a
generalization of a classical notion of a continuous family of maps between
locally compact spaces labeled by some other locally compact space.

This is based on the fact that for topological spaces $X,Y$ and $Z$ such that
$Z$ is Hausdorff and $X$ is locally compact (Hausdorff) we have
\[
\C(Z\times{X},Y)\approx\C\bigl(Z,\C(X,Y)\bigr),
\]
where ``$\approx$'' means homeomorphism and all spaces are taken with
compact-open topology \cite{jackson}. Thus a continuous family of maps
$X\to{Y}$ labeled by $Z$ can be represented by a continuous map
$Z\times{X}\to{Y}$.

\begin{defn}\label{DefQFM}
Let $\qs(A),\qs(B)$ and $\qs(C)$ be quantum spaces.
\begin{enumerate}
\item A \emph{quantum family of maps $\qs(C)\to\qs(B)$ labeled by $\qs(A)$}
is an element $\Psi\in\Mor(B,C\tens{A})$.
\item\label{DefQFM3} We say that $\Phi\in\Mor(B,C\tens{A})$ is the
\emph{quantum family
of all maps $\qs(C)\to\qs(B)$} if for any quantum space $\qs(D)$ and any
quantum family $\Psi$ of maps $\qs(C)\to\qs(B)$ labeled by $\qs(D)$
there exists a unique $\Lambda\in\Mor(A,D)$ such that the diagram
\[
\xymatrix{
B\ar[rr]^-{\Phi}\ar@{=}[d]&&C\tens{A}\ar[d]^{\id\tens\Lambda}\\
B\ar[rr]^-{\Psi}&&C\tens{D}}
\]
is commutative. In this case we say that $\qs(A)$ is the \emph{quantum space of
all maps $\qs(C)\to\qs(B)$.}
\item\label{DefQFM4} In the special case when $B=C$, we say that a quantum
family $\Psi\in\Mor(B,B\tens{A})$ is \emph{trivial} if $\Psi(b)=b\tens{I}$ for
all $b\in{B}$.
\end{enumerate}
The property of $(A,\Phi)$ described in \eqref{DefQFM3} will be referred to as
the \emph{universal property of $(A,\Phi)$.}
\end{defn}

\begin{rem}
\noindent\begin{enumerate}
\item Let $\Psi\in\Mor(B,C\tens{A})$ be a quantum family of maps
$\qs(C)\to\qs(B)$ labeled by $\qs(A)$. Assume that all three spaces are
classical, i.e.~$A,B$ and $C$ are commutative. Then $\qs(A)$ is a classical
locally compact space labeling a family of elements of
$\C\bigl(\qs(C),\qs(B)\bigr)$ and the map from $\qs(A)$ to
$\C\bigl(\qs(C),\qs(B)\bigr)$ is continuous for the compact-open topology.
\item It is clear that given two quantum spaces $\qs(B)$ and $\qs(C)$, the
quantum space of all maps $\qs(C)\to\qs(B)$ might not exists in the category
$\gQ$. This can happen even for classical spaces. Non existence of the
quantum space of all maps should be understood as meaning that this object
is not locally compact in the compact open-topology, rather than that it does
not exist at all.
\item If the quantum space $\qs(A)$ of all maps $\qs(C)\to\qs(B)$ exists and
$\Phi$ is the quantum family of all maps $\qs(C)\to\qs(B)$ then the pair
$(A,\Phi)$ is unique in the sense that if $(A',\Phi')$ is another such pair then
there exists an isomorphism $\Lambda\in\Mor(A,A')$ such that
$(\id\tens\Lambda)\comp\Phi=\Phi'$.
\end{enumerate}
\end{rem}

It was stated already in \cite{pseu} that the quantum space of all maps from
a finite quantum space (described by a finite dimensional $\cst$-algebra) to a
compact finite dimensional one (corresponding to a unital finitely generated 
$\cst$-algebra) always exists:

\begin{thm}\label{exSLW}
Let $\qs(B)$ and $\qs(C)$ be quantum spaces. Assume that $C$ is finite
dimensional and $B$ is finitely generated and unital. Then
\begin{enumerate}
\item\label{exSLW1} the quantum space $\qs(A)$ of all maps $\qs(C)\to\qs(B)$
exists.
\item\label{exSLW2} The $\cst$-algebra $A$ is unital and generated by
$\bigl\{(\omega\tens\id)\Phi(b)\st{b}\in{B},\:\omega\in{C^*}\bigr\}$, where
$\Phi\in\Mor(B,C\tens{A})$ is the quantum family of all maps $\qs(C)\to\qs(B)$.
\end{enumerate}
\end{thm}

\begin{proof}
Let $x_1,\ldots,x_N$ be generators of $B$. Since any element of $B$ is a linear
combinations of four unitaries, we can assume that $x_1,\ldots,x_N$ are
unitary. Let $(\mathcal{R}_t)_{t\in\mathcal{T}}$ be the complete list
of relations between $x_1,\ldots,x_N$, so that
\[
\bigl\langle{x_1,\ldots,x_N}\st\mathcal{R}_t(x_1,\ldots,x_N)=0,\:t\in\mathcal{T}
\bigr\rangle
\]
is a presentation of $B$. In particular some of the relations
$(\mathcal{R}_t)_{t\in\mathcal{T}}$ say that each $x_p$ is unitary.

The $\cst$-algebra $C$ can be written as a finite
direct sum of full matrix algebras:
\[
C=\bigoplus_{k=1}^KM_{n_k}.
\]
Let $\mathcal{A}$ be the $*$-algebra generated by elements
\[
\bigl\{{y_p^k}_{r,s}\st{p\in\{1,\ldots,N\}},\:
{k\in\{1,\ldots,K\}},\:r,s\in\{1,\ldots,n_k\}\bigr\}.
\]
with the relations
\[
\mathcal{R}_t
\left(
\begin{pmatrix}
Y_1^1\\
&\ddots\\
&&Y_1^K
\end{pmatrix}
,\ldots,
\begin{pmatrix}
Y_N^1\\
&\ddots\\
&&Y_N^K
\end{pmatrix}
\right)=0,\qquad(t\in\mathcal{T}),
\]
where $Y_p^k$ is the matrix $\bigl({y_p^k}_{r,s}\bigr)_{r,s=1,\ldots,n_k}$.

The relation saying that $x_p$ is unitary guarantees that
\[
\sum_{q=1}^{n_k}({y_p^k}_{q,r})^*({y_p^k}_{q,r})=\delta_{r,s}I_{\mathcal{A}}.
\]
In particular, for any Hilbert space representation $\pi$ of $\mathcal{A}$ the
norm $\bigl\|\pi({y_p^k}_{r,s})\bigr\|\leq{1}$.

This implies that for any $c\in\mathcal{A}$ the quantity
\begin{equation}\label{supH}
\sup_{\pi}\bigl\|\pi(c)\bigr\|
\end{equation}
(where the supremum is taken over all Hilbert space representations of
$\mathcal{A}$) is finite. Now standard procedure leads to construction of the
universal enveloping $\cst$-algebra $A$ of $\mathcal{A}$.

Let us define a map $\Phi:B\to{C\tens{A}}$ as sending $x_p$ to the image in
$C\tens{A}$ of the matrix
\[
\begin{pmatrix}
Y_p^1\\
&\ddots\\
&&Y_p^K
\end{pmatrix}\in\C\tens\mathcal{A}.
\]
After a moment of reflection, we see that $(A,\Phi)$ has the required universal
property, so that $\qs(A)$ is the space of all maps $\qs(C)\to\qs(B)$ and
$\Phi\in\Mor(B,C\tens{A})$ is the quantum family of all these maps.

The second part of Statement \eqref{exSLW2} follows from the obvious
observation that
\[
\bigl\{(\omega\tens\id)\Phi(b)\st{b}\in{B},\:\omega\in{C^*}\bigr\}
\]
contains the images of the generators of $\mathcal{A}$ in $A$. Alternatively,
one can prove this using the uniqueness of $(A,\Phi)$.
\end{proof}

\begin{rem}
\noindent\begin{enumerate}
\item In case when $C$ is the algebra $M_n$ of $n\times{n}$ complex matrices,
the $\cst$-algebra $A$ defined in Theorem \ref{exSLW} coincides with
$W_n(M_n)$, where $W_n$ is the left adjoint of the functor $D\mapsto{M_n(D)}$
(cf.~\cite[Page 174]{invlc}). Note, however that our morphism sets are
different that those used in \cite{invlc}.
\item It was noticed in \cite{pseu} that the quantum space of maps between
quantum spaces might be very interesting even if the two quantum spaces are
finite classical spaces. For example if we take $B=C=\CC^2$ then $\qs(B)=\qs(C)$
is the two-point space. The family of all maps $\qs(C)\to\qs(B)$ has
four elements. However the corresponding \emph{quantum} family is infinite in
the sense that the corresponding $\cst$-algebra is infinite dimensional. It is
amusing to check that in this case it is isomorphic to the $\cst$-algebra of
all continuous functions from an interval to $M_2$ whose values at the
endpoints are diagonal (\cite{pseu}, \cite[Section 2.$\beta$]{ncg}).
\end{enumerate}
\end{rem}

We shall now introduce the notion of composition of quantum families of
maps. Let $A_1$, $A_2$, $B$, $C$ and $D$ be $\cst$-algebras and let
\[
\Psi_1\in\Mor(C,D\tens{A_1}),\qquad\Psi_2\in\Mor(B,C\tens{A_2})
\]
be quantum families of maps $\qs(D)\to\qs(C)$ and $\qs(C)\to\qs(B)$ labeled by
$\qs(A_1)$ and $\qs(A_2)$ respectively. We define the quantum family
$\Psi_1\vt\Psi_2$ of maps $\qs(D)\to\qs(B)$ labeled by $\qs(A_1\tens{A_2})$ by
\[
\Psi_1\vt\Psi_2=(\Psi_1\tens\id)\comp\Psi_2.
\]
This quantum family of maps will be called the \emph{composition} of the
families $\Psi_1$ and $\Psi_2$. We shall also refer to the operation taking
$\Psi_1$ and $\Psi_2$ to $\Psi_1\vt\Psi_2$ as the \emph{operation of
composition} of quantum families of maps. In case the families are classical,
i.e.~the $\cst$-algebras $A_1$ and $A_2$ are commutative, the family
$\Psi_1\vt\Psi_2$ is a classical family consisting of all compositions of
members of $\Psi_1$ and $\Psi_2$.

The crucial property of composition of quantum families of maps is that it is
associative:

\begin{prop}\label{ASC}
Let $A_1$, $A_2$, $A_3$, $B_1$, $B_2$, $C$ and $D$ be $\cst$-algebras and let
\[
\begin{split}
\Psi_1&\in\Mor(B_1,C\tens{A_1}),\\
\Psi_2&\in\Mor(B_2,B_1\tens{A_2}),\\
\Psi_3&\in\Mor(D,B_2\tens{A_3})
\end{split}
\]
be quantum families of maps. Then
\[
\Psi_1\vt(\Psi_2\vt\Psi_3)=(\Psi_1\vt\Psi_2)\vt\Psi_3.
\]
\end{prop}

\begin{proof}
This is a simple computation:
\[
\begin{split}
\Psi_1\vt(\Psi_2\vt\Psi_3)
&=(\Psi_1\tens\id)\comp(\Psi_2\vt\Psi_3)\\
&=(\Psi_1\tens\id\tens\id)\comp(\Psi_2\tens\id)\comp\Psi_3\\
&=\bigl(\bigl[(\Psi_1\tens\id)\comp\Psi_2\bigr]\tens\id\bigr)\comp\Psi_3\\
&=\bigl([\Psi_1\vt\Psi_2]\tens\id\bigr)\comp\Psi_3\\
&=(\Psi_1\vt\Psi_2)\vt\Psi_3.
\end{split}
\]
\end{proof}

\section{Quantum semigroup structure}\label{qss}

In this section we shall analyze the structure of the quantum space of all
maps from a finite quantum space to itself. Thus let $M$ be a finite
dimensional $\cst$-algebra. Then Theorem \ref{exSLW} guarantees that there
exists the quantum space of all maps $\qs(M)\to\qs(M)$. Let us denote the
corresponding $\cst$-algebra by $A$ and let $\Phi\in\Mor(M,M\tens{A})$ be
the quantum family of all maps $\qs(M)\to\qs(M)$. This notation will be kept
throughout this section.

The universal property of $(A,\Phi)$ will provide us with rich structure on
$A$ or, more appropriately, on $\qs(A)$.

\begin{thm}\label{qstr}
\noindent\begin{enumerate}
\item\label{qstr1} There exists a unique morphism $\Delta\in\Mor(A,A\tens{A})$
such that
\begin{equation}\label{PhiDel}
(\Phi\tens\id)\comp\Phi=(\id\tens\Delta)\comp\Phi.
\end{equation}
\item\label{qstr2} The morphism $\Delta$ satisfies
\begin{equation}\label{coas}
(\Delta\tens\id)\comp\Delta=(\id\tens\Delta)\comp\Delta.
\end{equation}
\item\label{qstr3} There exists a unique character $\epsilon$ of $A$ such that
\begin{equation}\label{PhiEps}
(\id\tens\epsilon)\comp\Phi=\id.
\end{equation}
\item\label{qstr4} The character $\epsilon$ satisfies
\begin{equation}\label{coun}
(\id\tens\epsilon)\comp\Delta=(\epsilon\tens\id)\comp\Delta=\id.
\end{equation}
\end{enumerate}
\end{thm}

\begin{proof}
Let us consider the quantum family $\Phi\vt\Phi$ of maps $\qs(M)\to\qs(M)$. It
is labeled by $\qs(A\tens{A})$ and the universal property of $(A,\Phi)$
implies that there exists a unique $\Delta\in\Mor(A,A\tens{A})$ such that
\[
(\id\tens\Delta)\comp\Phi=\Phi\vt\Phi.
\]
This is precisely \eqref{PhiDel}.

To prove \eqref{coas} we use \eqref{PhiDel} to compute
$(\Phi\tens\id\tens\id)\comp(\Phi\tens\id)\comp\Phi$ in two ways:
\[
\begin{split}
(\Phi\tens\id\tens\id)\comp(\Phi\tens\id)\comp\Phi
&=(\Phi\tens\id\tens\id)\comp(\id\tens\Delta)\comp\Phi\\
&=(\id\tens\id\tens\Delta)\comp(\Phi\tens\id)\comp\Phi\\
&=(\id\tens\id\tens\Delta)\comp(\id\tens\Delta)\comp\Phi\\
\end{split}
\]
and
\[
\begin{split}
(\Phi\tens\id\tens\id)\comp(\Phi\tens\id)\comp\Phi
&=\Bigl(\bigl[(\Phi\tens\id)\comp\Phi\bigr]\tens\id\Bigr)\comp\Phi\\
&=\Bigl(\bigl[(\id\tens\Delta)\comp\Phi\bigr]\tens\id\Bigr)\comp\Phi\\
&=(\id\tens\Delta\tens\id)\comp(\Phi\tens\id)\comp\Phi\\
&=(\id\tens\Delta\tens\id)\comp(\id\tens\Delta)\comp\Phi.
\end{split}
\]
Let $\omega$ be a functional on $M$. Applying $(\omega\tens\id\tens\id)$ to
both sides of the equation
\[
(\id\tens\id\tens\Delta)\comp(\id\tens\Delta)\comp\Phi
=(\id\tens\Delta\tens\id)\comp(\id\tens\Delta)\comp\Phi
\]
we obtain
\[
\bigr[(\id\tens\Delta)\comp\Delta\bigl]
\bigl((\omega\tens\id)\Phi(m)\bigr)=
\bigr[(\Delta\tens\id)\comp\Delta\bigl]
\bigl((\omega\tens\id)\Phi(m)\bigr).
\]
for any $m\in{M}$. Thus formula \eqref{coas} follows form Proposition
\ref{exSLW}\eqref{exSLW2}.

Statement \eqref{qstr3} follows again from the universal property of
$(A,\Phi)$ (Definition \ref{DefQFM}\eqref{DefQFM3} with $B=C=M$). More
precisely we take $D=\CC$ and canonically identify $M$ with
$M\tens{D}$. Then we take $\Psi$ to be the identity morphism $M\to{M\tens{D}}$.
The universal property of $(A,\Phi)$ guarantees that there exists a unique
$\epsilon\in\Mor(A,D)$ satisfying \eqref{PhiEps}.

The proof of Statement \eqref{qstr4} is similar to that of Statement
\eqref{qstr2}. Using \eqref{PhiDel} and \eqref{PhiEps} we arrive at
\[
\Bigl(\id\tens\bigl[(\epsilon\tens\id)\comp\Delta\bigr]\Bigr)\comp\Phi
=\Phi=
\Bigl(\id\tens\bigl[(\id\tens\epsilon)\comp\Delta\bigr]\Bigr)\comp\Phi.
\]
Applying $(\omega\tens\id\tens\id)$ to both sides, we obtain
\[
\bigl[(\epsilon\tens\id)\comp\Delta\bigr]
\bigl((\omega\tens\id)\Phi(m)\bigr)
=(\omega\tens\id)\Phi(m)=
\bigl[(\id\tens\epsilon)\comp\Delta\bigr]
\bigl((\omega\tens\id)\Phi(m)\bigr)
\]
for any $m\in{M}$ which proves \eqref{coun}.
\end{proof}

The morphisms $\Delta\in\Mor(A,A\tens{A})$ and $\epsilon\in\Mor(A,\CC)$ are
called the \emph{comultiplication} and the \emph{counit} of $A$ respectively.
They endow $A$ with the structure of a \emph{quantum semigroup} with unit as
defined below.

\begin{defn}
\noindent\begin{enumerate}
\item A pair $(B,\Delta_B)$ consisting of a $\cst$-agebra and a morphism
$\Delta_B\in\Mor(B,B\tens{B})$ is called a \emph{quantum semigroup} if
$\Delta_B$ is \emph{coassociative}, i.e.
\[
(\Delta_B\tens\id)\comp\Delta_B=(\id\tens\Delta_B)\comp\Delta_B.
\]
\item A quantum semigroup $(B,\Delta_B)$ has a unit if $B$ admits a character
$\epsilon_B$ satisfying
\[
(\epsilon_B\tens\id)\comp\Delta_B=\id=(\id\tens\epsilon_B)\comp\Delta_B.
\]
\item Let $(B,\Delta_B)$ and $(C,\Delta_C)$ be quantum semigroups. An element
$\Lambda\in\Mor(B,C)$ is a \emph{morphism of quantum semigroups} (or a
\emph{quantum semigroup morphism}) if it satisfies
\[
(\Lambda\tens\Lambda)\comp\Delta_B=\Delta_C\comp\Lambda.
\]
\item An action of a quantum semigroup $(B,\Delta_B)$ on a quantum space
$\qs(C)$ is a morphism $\Psi\in\Mor(C,C\tens{B})$ satisfying
\[
(\Psi\tens\id)\comp\Psi=(\id\tens\Delta_B)\comp\Psi.
\]
\end{enumerate}
\end{defn}

We shall denote the quantum semigroup $(A,\Delta)$ constructed in Theorem
\ref{qstr} by $\QMap\bigl(\qs(M)\bigr)$. The
quantum family $\Phi\in\Mor(M,M\tens{A})$ of all maps $\qs(M)\to\qs(M)$ is then
an action of $\QMap\bigl(\qs(M)\bigr)$ on
the quantum space $\qs(M)$. Since $A$ is a unital $\cst$-algebra (Theorem
\ref{exSLW}\eqref{exSLW2}), the quantum semigroup $\QMap\bigl(\qs(M)\bigr)$ is
compact.

\begin{rem}
\noindent\begin{enumerate}
\item The coassociativity of $\Delta$ as derived in the proof of Theorem
\ref{qstr} is, in fact, a consequence of the associativity of composition of
quantum families of maps (Proposition \ref{ASC}).
\item The semigroup structure on $\QMap\bigl(\qs(M)\bigr)$ and the operation of
composition of quantum families of maps are related in a very natural way. Let
$B$ and $C$ be $\cst$-algebras and let
\[
\Psi_B\in\Mor(M,M\tens{B}),\qquad\Psi_C\in\Mor(M,M\tens{C})
\]
be quantum families of maps $\qs(M)\to\qs(M)$ labeled by $\qs(B)$ and $\qs(C)$
respectively. Finally let $\Lambda_B\in\Mor(A,B)$ and $\Lambda_C\in\Mor(A,C)$
be the unique morphisms satisfying
\[
(\id\tens\Lambda_B)\comp\Phi=\Psi_B,\qquad(\id\tens\Lambda_C)\comp\Phi=\Psi_C.
\]
Then
\[
\Psi_B\vt\Psi_C=(\id\tens\Lambda_B\tens\Lambda_C)\comp(\id\tens\Delta)\comp\Phi.
\]
Clearly the morphisms $\Lambda_B$ and $\Lambda_C$ describe the ``inclusions''
of the two considered families into the quantum family of all maps, and their
composition is the composition in the semigroup of all maps.
\end{enumerate}
\end{rem}

The coassociativity of $\Delta$ means, in particular that the operation of
\emph{convolution product} of continuous functionals on $A$, defined by
\[
\phi*\psi=(\phi\tens\psi)\comp\Delta
\]
for $\phi,\psi\in{A^*}$, is associative. Note that the counit $\epsilon$
is, by Theorem \ref{qstr}\eqref{qstr4}, a neutral element for the
convolution product. The notion of convolution product enters the
formulation of the next theorem.

Also in the next theorem we shall use the concept of natural topology on the
set of morphisms between $\cst$-algebras defined as follows: let $B$ and $C$
be $\cst$-algebras. The \emph{natural topology on $\Mor(B,C)$} is the weakest
topology such that for any $b\in{B}$ the maps
\[
\Mor(B,C)\ni\Psi\longmapsto\Psi(b)\in\M(C)
\]
is strictly continuous (\cite[p.~491]{gen}). If $B=\Cinf(Y)$ and $C=\Cinf(X)$
are commutative, the set $\Mor(B,C)$ can be identified with $\C(X,Y)$ and the
natural topology is the compact-open topology. In the relevant case, when
$B=C=M$ is a finite dimensional $\cst$-algebra, this topology is the one
inherited by $\Mor(M,M)$ form the space of all linear maps $M\to{M}$ which is
a finite dimensional vector space.

\begin{thm}\label{bi}
\noindent\begin{enumerate}
\item\label{bi1} There exists a canonical bijection $\Theta$ between the
space of all characters of $A$ and the set $\Mor(M,M)$ given by
\[
\xymatrix{
M\ar[rr]^-{\Phi}\ar@{=}[d]&&M\tens{A}\ar[d]^{\id\tens\lambda}\\
M\ar[rr]^-{\Theta(\lambda)}&&M}
\]
\item\label{bi2} The operation of convolution endows the set of all characters
of $A$ with the structure of a unital semigroup and $\Theta$ is an isomorphism
of semigroups.
\item\label{bi3} $\Theta$ is a homeomorphism for the weak$^*$ topology on the
space of characters of $A$ and the natural topology on $\Mor(M,M)$.
\item\label{bi4} $\Theta$ maps the set of all convolution invertible characters
onto the group of all automorphisms of $M$ and is an isomorphism of topological
groups.
\end{enumerate}
\end{thm}

\begin{proof}
Statement \eqref{bi1} follows from the universal property of $(A,\Phi)$: any
morphism $\Lambda\in\Mor(M,M)$ gives a character $\lambda$ on $A$ such
that $(\id\tens\lambda)\comp\Phi=\Lambda$. Conversely, any character $\lambda$
on $A$ defines $\Lambda=(\id\tens\lambda)\comp\Phi\in\Mor(M,M)$. These
correspondences are inverse to one another by the universal property of
$(A,\Phi)$.

The fact that $\Theta$ is a semigroup isomorphism follows from inspection of
the following commutative diagram:
\[
\xymatrix{
M\ar[rr]^-{\Phi}\ar@{=}[d]&&M\tens{A}\ar[rr]^{\id\tens\Delta}
&&M\tens{A}\tens{A}\ar@{=}[d]
\ar@<5ex>@/^5ex/[ddd]^{\id\tens\lambda\tens\mu}\\
M\ar[rr]^-{\Phi}\ar@{=}[d]&&M\tens{A}\ar[d]^{\id\tens\mu}\ar[rr]^{\Phi\tens\id}
&&M\tens{A}\tens{A}\ar[d]^{\id\tens\id\tens\mu}\\
M\ar[rr]^-{\Theta(\mu)}&&M\ar[rr]^-{\Phi}\ar@{=}[d]
&&M\tens{A}\ar[d]^{\id\tens\lambda}\\
&&M\ar[rr]^-{\Theta(\lambda)}&&M
}
\]
We obtain the equality
$\bigl(\id\tens[\lambda*\mu]\bigr)\comp\Phi=\Theta(\lambda)\comp\Theta(\mu)$,
so that
\[
\Theta(\lambda*\mu)=\Theta(\lambda)\comp\Theta(\mu).
\]

To prove Statement \eqref{bi3} it is enough to show that $\Theta$ is continuous
for the weak$^*$ topology on $A$ and the natural topology on $\Mor(M,M)$. The
fact that $\Theta$ is a homeomorphism will follow, because both spaces are
compact and $\Theta$ is a bijection.

Take $m\in{M}$ and let $(\lambda_\alpha)$ be a net of characters of $A$,
weak$^*$ convergent to $\lambda$ (which must then be a character of $A$). It is
clear that $(\id\tens\lambda_\alpha)\Phi(m)$ is norm convergent to
$(\id\tens\lambda)\Phi(m)$ ($M$ is finite dimensional). This means that
$\bigl(\Theta(\lambda_\alpha)\bigr)(m)$ depends continuously on $\alpha$. Since
$m$ is arbitrary, this means that $\Theta(\lambda_\alpha)$ varies continuously
for the natural topology on $\Mor(M,M)$.

Statement \eqref{bi3} is a consequence of \eqref{bi2} and \eqref{bi3}.
\end{proof}

The next statement says, in particular, that the quantum semigroup
$\QMap\bigl(\qs(M)\bigr)$ is not to small, because it always contains the
semigroup $\Mor(M,M)$. Containment must be understood in the sense of
noncommutative geometry.

\begin{cor}
The Gelfand transform of $A$ is an epimorphism onto $\C\bigl(\Mor(M,M)\bigr)$
and is a morphism of quantum semigroups.
\end{cor}

An interesting question is whether the quantum semigroup
$\QMap\bigl(\qs(M)\bigr)$ is significantly bigger than the semigroup
$\Mor(M,M)$. A partial answer to this question is given by the following
proposition.

\begin{prop}\label{ergodic}
The action $\Phi$ of $\QMap\bigl(\qs(M)\bigr)$ on $\qs(M)$ is ergodic: for any
$m\in{M}$, the condition that $\Phi(m)=m\tens{I}$ implies $m\in\CC{I}$.
\end{prop}

\begin{proof}
Consider the quantum family of maps $\qs(M)\to\qs(M)$ labeled by $\qs(M)$
given by $\Psi\in\Mor(M,M\tens{M})$,
\[
\Psi(m)=I\tens{m}.
\]
By the universal property of $(A,\Phi)$ there is an element
$\Lambda\in\Mor(A,M)$ such that $I\tens{m}=(\id\tens\Lambda)\Phi(m)$. Now if
$\Phi(m)=m\tens{I}$ then $I\tens{m}=(\id\tens\Lambda)(m\tens{I})=m\tens{I}$ and
it follows that $m\in\CC{I}$.
\end{proof}

Note that the semigroup $\Mor(M,M)$ need not be ergodic in the sense that there
can be nontrivial elements of $M$ fixed under every morphism.

\begin{prop}\label{qsm}
Let $B$ be a $\cst$-algebra and let $\Psi\in\Mor(M,M\tens{B})$ be a quantum
family of maps $\qs(M)\to\qs(M)$ labeled by $\qs(B)$. Assume that there exists a
morphism $\Delta_B\in\Mor(B,B\tens{B})$ such that
\[
(\id\tens\Delta_B)\comp\Psi=(\Psi\tens\id)\comp\Psi
\]
and let $\Lambda\in\Mor(A,B)$ be the unique morphism such that
$(\id\tens\Lambda)\comp\Phi=\Psi$. Then $\Lambda$ satisfies
\[
(\Lambda\tens\Lambda)\comp\Delta=\Delta_B\comp\Lambda.
\]
\end{prop}

\begin{proof}
Looking at the commutative diagram
\[
\xymatrix{
M\ar[rr]^-{\Phi}\ar@{=}[d]&&
M\tens{A}\ar[rr]^-{\id\tens\Delta}&&
M\tens{A}\tens{A}\ar@{=}[d]\\
M\ar[rr]^-{\Phi}\ar@{=}[d]&&
M\tens{A}\ar[rr]^-{\Phi\tens\id}
\ar[d]^{\id\tens\Lambda}&&
M\tens{A}\tens{A}\ar[d]^{\id\tens\Lambda\tens\Lambda}\\
M\ar[rr]^-{\Psi}&&M\tens{B}\ar[rr]^-{\Psi\tens\id}&&M\tens{B}\tens{B}
}
\]
we find that $(\id\tens\Lambda\tens\Lambda)\comp(\id\tens\Delta)\comp\Phi
=(\Psi\tens\id)\comp(\id\tens\Lambda)\comp\Phi$. But
$(\id\tens\Lambda)\comp\Phi=\Psi$ and using this fact twice we obtain
\[
\begin{split}
(\id\tens\Lambda\tens\Lambda)\comp(\id\tens\Delta)\comp\Phi
&=(\Psi\tens\id)\comp(\id\tens\Lambda)\comp\Phi\\
&=(\Psi\tens\id)\comp\Psi\\
&=(\id\tens\Delta_B)\comp\Psi\\
&=(\id\tens\Delta_B)\comp(\id\tens\Lambda)\comp\Phi.
\end{split}
\]
Take $\omega\in{M^*}$ and let us apply $\omega\tens\id\tens\id$ to both sides
of the above equality. We find that for any $m\in{M}$
\[
(\Delta_B\comp\Lambda)\bigl((\omega\tens\id)\Phi(m)\bigr)
=\bigl((\Lambda\tens\Lambda)\comp\Delta\bigr)
\bigl((\omega\tens\id)\Phi(m)\bigr)
\]
and our result follows from Theorem \ref{exSLW}\eqref{exSLW2}.
\end{proof}

\section{Invariant states}\label{invS}

We shall retain the notation introduced in Section \ref{qss}. Thus $M$ is a
finite dimensional $\cst$-algebra and $\QMap\bigl(\qs(M)\bigr)=(A,\Delta)$ is
the quantum semigroup of all maps $\qs(M)\to\qs(M)$. The action of
$\QMap\bigl(\qs(M)\bigr)$ on $\qs(M)$ will, as before, be denoted by $\Phi$.

\begin{defn}
Let $\Psi\in\Mor(M,M\tens{B})$ be a quantum family of maps $\qs(M)\to\qs(M)$
labeled by $\qs(B)$ and let $\omega$ be a state on $M$. We say that $\omega$ is
\emph{invariant} for $\Psi$ if
\[
(\omega\tens\id)\Psi(m)=\omega(m)I
\]
for all $m\in{M}$. We also say that the quantum family of maps $\Psi$
\emph{preserves} the state $\omega$.
\end{defn}

In the next proposition we shall use the notion of composition of quantum
families of maps defined in Section \ref{qspm}.

\begin{prop}\label{compw}
Let $B$ and $C$ be $\cst$-algebras and let $\Psi_B\in\Mor(M,M\tens{B})$,
$\Psi_C\in\Mor(M,M\tens{C})$ be quantum families of maps $\qs(M)\to\qs(M)$
labeled by $\qs(B)$ and $\qs(C)$ respectively. Let $\omega$ be a state on $M$
which is invariant for both $\Psi_B$ and $\Psi_C$. Then $\omega$ is invariant
for the composition $\Psi_B\vt\Psi_C$.
\end{prop}

\begin{proof}
The proof of this result is purely computational. To make the computations
more transparent let us denote them maps
\[
\begin{split}
\CC\ni{z}\longmapsto{z}I&\in{B},\\
\CC\ni{z}\longmapsto{z}I&\in{C},\\
\CC\ni{z}\longmapsto{z}I&\in{B\tens{C}},\\
\end{split}
\]
by $\Eta_B$, $\Eta_C$ and $\Eta_{B\tens{C}}$ respectively. Then we have
\[
(\omega\tens\id)\comp\Psi_B=\Eta_B\comp\omega,\qquad
(\omega\tens\id)\comp\Psi_C=\Eta_C\comp\omega
\]
and $\Eta_B\tens\Eta_C=\Eta_{B\tens{C}}$ (with the identification
$\CC\tens\CC=\CC$).

Now
\[
\begin{split}
(\omega\tens\id)\comp(\Psi_B\vt\Psi_C)
&=(\omega\tens\id\tens\id)\comp(\Psi_B\tens\id)\comp\Psi_C\\
&=\bigl(\bigl[(\omega\tens\id)\comp\Psi_B\bigr]\tens\id\bigr)\comp\Psi_C\\
&=\bigl([\Eta_B\comp\omega]\tens\id\bigr)\comp\Psi_C\\
&=(\Eta_B\tens\id)\comp(\omega\tens\id)\comp\Psi_C\\
&=(\Eta_B\tens\id)\comp\Eta_C\comp\omega\\
&=(\Eta_B\tens\Eta_B)\comp\omega
=\Eta_{B\tens{C}}\comp\omega,
\end{split}
\]
so that for any $m\in{M}$ we have
$(\omega\tens\id)(\Psi_B\vt\Psi_C)(m)=\omega(m)I$.
\end{proof}

One can ask if there are any states on $M$ preserved by the action of
$\QMap\bigl(\qs(M)\bigr)$, i.e.~if there is a state $\omega$ on $M$ such that
for any $m\in{M}$ we have
\begin{equation}\label{invar}
(\omega\tens\id)\Phi(m)=\omega(m)I.
\end{equation}
It turns out that, unless $M$ is one dimensional, there are no such states.
This result is proved in the same way as Proposition \ref{ergodic}.

\begin{prop}\label{erg-omega}
Let $\omega$ be a state on $M$. Assume that $\omega$ is invariant under the
action of $\QMap\bigl(\qs(M)\bigr)$. Then $M$ is one dimensional.
\end{prop}

\begin{proof}
Let $\Lambda\in\Mor(A,M)$ be such that for
any $m\in{M}$ we have $(\id\tens\Lambda)\Phi(m)=I\tens{m}$ (see proof of
Proposition \ref{ergodic} for the existence of $\Lambda$). Applying $\Lambda$
to both sides of \eqref{invar} we obtain
\[
\Lambda\bigl((\omega\tens\id)\Phi(m)\bigr)=\omega(m)\Lambda(I).
\]
The right hand side of this equality is $\omega(m)I$, while the left hand side
is
\[
(\omega\tens\id)\bigl((\id\tens\Lambda)\Phi(m)\bigr)
=(\omega\tens\id)(I\tens{m})=m,
\]
so that $m=\omega(m)I$ for any $m\in{M}$.
\end{proof}

In the next theorem we shall describe the quantum subsemigroups of
$\QMap\bigl(\qs(M)\bigr)$ preserving a given state $\omega$ on $M$. The
strategy is to define the smallest ideal that must be included in the kernel of
any morphism from $A$ to any $\cst$-algebra such that the resulting quantum
family of maps $\qs(M)\to\qs(M)$ preserves $\omega$. Then the quotient of $A$
by this ideal gives a quantum family of maps $\qs(M)\to\qs(M)$ which is
universal for all quantum families preserving $\omega$. This universality will
give a comultiplication on the quotient $\cst$-algebra.

\begin{thm}\label{Qw}
Let $\omega$ be a state on $M$ and let $J$ be the ideal generated by the set
\[
\Bigl\{(\omega\tens\id)\Phi(m)-\omega(m)I\st{m\in{M}}\Bigr\}.
\]
Let $\dot{A}$ be the quotient $A/J$, let $\pi:A\to\dot{A}$ be the canonical
epimorphism and let $\dot{\Phi}=(\id\tens\pi)\comp\Phi$. Then
\begin{enumerate}
\item\label{invw} the state $\omega$ is invariant for the quantum
family $\dot{\Phi}\in\Mor\bigl(M,M\tens\dot{A}\bigr)$ of maps $\qs(M)\to\qs(M)$
labeled by $\qs\bigl(\dot{A}\bigr)$.
\item\label{univw} For any $\cst$-algebra $B$ and any quantum family
$\Psi\in\Mor(M,M\tens{B})$ of maps $\qs(M)\to\qs(M)$ labeled by $\qs(B)$ such
that $\omega$ is invariant for $\Psi$ there exists a unique morphism
$\Lambda\in\Mor\bigl(\dot{A},B\bigr)$
such that $\Psi=(\id\tens\Lambda)\comp\dot{\Phi}$.
\item\label{exdDel} There exists a unique
$\dot{\Delta}\in\Mor\bigl(\dot{A},\dot{A}\tens\dot{A}\bigr)$
such that
\begin{equation}\label{dPhi}
\bigl(\dot{\Phi}\tens\id\bigr)\comp\dot{\Phi}=\bigl(\id\tens\dot{\Delta}
\bigr)\comp\dot{\Phi}.
\end{equation}
\item The morphism $\dot{\Delta}$ is coassociative and
$\bigl(\dot{A},\dot{\Delta}\bigr)$ is a compact quantum semigroup with unit.
$\dot{\Phi}$ is an action of $\bigl(\dot{A},\dot{\Delta}\bigr)$ on $\qs(M)$.
\item\label{piqsm} $\pi$ is a quantum semigroup morphism.
\item\label{Lqsm} For any quantum semigroup $(B,\Delta_B)$ acting on $\qs(M)$
with action $\Phi_B\in\Mor(M,M\tens{B})$ preserving $\omega$, the unique
morphism
$\Lambda\in\Mor\bigl(\dot{A},B\bigr)$ such that
$\Phi_B=(\id\tens\Lambda)\comp\dot{\Phi}$ is a quantum semigroup morphism.
\end{enumerate}
\end{thm}

\begin{proof}
Statement \eqref{invw} is almost obvious. Since $\pi$ sends each element of the
form
\begin{equation}\label{form}
(\omega\tens\id)\Phi(m)-\omega(m)I
\end{equation}
to $0$, we see that we have
\[
(\omega\tens\id)\dot{\Phi}(m)=\omega(m)I
\]
for all $m\in{M}$.

Let $B$ be a $\cst$-algebra and let $\Psi\in\Mor(M,M\tens{B})$ be a quantum
family of maps $\qs(M)\to\qs(M)$ preserving $\omega$. Then there is a unique
map $\Lambda_0\in\Mor(A,B)$ such that $(\id\tens\Lambda_0)\comp\Phi=\Psi$.
Moreover all elements of the form \eqref{form} are mapped to $0$ by
$\Lambda_0$. Therefore $J\subset\ker{\Lambda_0}$. This guarantees that there
is a unique $\Lambda\in\Mor\bigl(\dot{A},B\bigr)$ such that
$\Lambda\comp\pi=\Lambda_0$. This means that
\[
(\id\tens\Lambda)\comp\dot{\Phi}
=(\id\tens\Lambda)\comp(\id\tens\pi)\comp\Phi
=(\id\tens\Lambda_0)\comp\Phi=\Psi
\]
and Statement \eqref{univw} is proven.

To prove statement \eqref{exdDel} note that by Proposition \ref{compw} the
composition $\dot{\Phi}\vt\dot{\Phi}$ is a quantum family of maps
$\qs(M)\to\qs(M)$ labeled by $\qs\bigl(\dot{A}\tens\dot{A}\bigr)$
for which $\omega$ is invariant. Therefore there exists a unique
$\dot{\Delta}\in\Mor\bigl(\dot{A},\dot{A}\tens\dot{A}\bigr)$ such that
\begin{equation}\label{unidD}
\bigl(\id\tens\dot{\Delta}\bigr)\comp\dot{\Phi}=\dot{\Phi}\vt\dot{\Phi}.
\end{equation}
By definition of $\dot{\Phi}\vt\dot{\Phi}$ we have
\[
\bigl(\id\tens\dot{\Delta}\bigr)\comp\dot{\Phi}
=\bigl(\dot{\Phi}\tens\id\bigr)\comp\dot{\Phi}.
\]
In fact \eqref{dPhi} is equivalent to \eqref{unidD}, so that $\dot{\Delta}$ is
unique.

The proof that $\bigl(\dot{A},\dot{\Delta}\bigr)$ is a quantum semigroup with
unit can be copied verbatim from that of Theorem \ref{qstr} and supplying
$\Phi$'s and $\Delta$'s with dots. The necessary facts for this are:
\begin{itemize}
\item the uniqueness of $\Lambda$ in Statement \eqref{univw},
\item $\dot{A}$ is generated by
$\bigl\{(\eta\tens\id)\dot{\Phi}(m)\st{m\in{M}},\:\eta\in{M^*}\bigr\}$ (this is
clearly seen from the analogous property for $A$ and the fact that $\pi$ is an
epimorphism).
\end{itemize}
The fact that $\dot{\Phi}$ is an action is obvious in view of \eqref{dPhi}.

Statement \eqref{piqsm} follows from Proposition \ref{qsm}. The proof of
\eqref{Lqsm} is analogous to the proof of Proposition \ref{qsm}. Again we need
to use uniqueness of $\Lambda$ stated in \eqref{univw}.
\end{proof}

The quantum semigroup $\bigl(\dot{A},\dot{\Delta}\bigr)$ constructed in Theorem
\ref{Qw} will be denoted by the symbol $\QMap^{\omega}\bigl(\qs(M)\bigr)$. This
is the quantum semigroup of all maps $\qs(M)\to\qs(M)$ preserving $\omega$.

\begin{rem}\label{ost}
\noindent\begin{enumerate}
\item In Theorem \ref{Qw} we obtained the comultiplication on $\dot{A}$ through
its universal property. However it is possible to show directly that
\begin{equation}\label{JAAJ}
\Delta(J)\subset{J\tens{A}}+A\tens{J},
\end{equation}
which gives a unique $\dot{\Delta}$ on
$\dot{A}$ such that $\pi$ is a quantum semigroup morphism. This can be seen for
example by choosing a basis $(m_k)_{k=1,\ldots,N}$ of $M$ and denoting by
$(a_{kl})_{k,l=1,\ldots,N}$ the elements of $A$ such that
$\Phi(m_l)=\sum\limits_{k=1}^Nm_k\tens{a_{kl}}$. Then the ideal $J$ is
generated by the elements
$X_l=\sum\limits_{k=1}^N\omega(m_k)a_{kl}-\omega(m_l)I$. One finds that
$\Delta(X_l)=\sum\limits_{p=1}^NX_p\tens{a_{pl}}+I\tens{X_l}$ and this suffices
for \eqref{JAAJ}.
\item Since $J\subset\ker{\epsilon}$, the counit of
$\bigl(\dot{A},\dot{\Delta}\bigr)$ composed with $\pi$ is $\epsilon$.
\item\label{ost3} From the fact that $\pi$ is an epimorphism and Theorem \ref{exSLW}
\eqref{exSLW2} it follows that the set
\[
\bigl\{(\eta\tens\id)\dot{\Phi}(m)\st{m\in{M}},\:\eta\in{M^*}\bigr\}
\]
generates $\dot{A}$ as a $\cst$-algebra.
\item The construction of $\QMap^\omega\bigl(\qs(M)\bigr)$ can be performed for
any $\omega\in{M^*}$, not necessarily a state. Theorem \ref{Qw} remains true in
that situation.
\end{enumerate}
\end{rem}

\section{Quantum commutants}\label{qc}

In this section we shall use the notation introduced in Section \ref{qss}. Let
us consider two $\cst$-algebras $B$ and $C$ and two quantum families
\[
\Psi_B\in\Mor(M,M\tens{B}),\qquad\Psi_C\in\Mor(M,M\tens{C})
\]
of maps $\qs(M)\to\qs(M)$ labeled by $\qs(B)$ and $\qs(C)$ respectively. We
shall say that $\Psi_B$ \emph{commutes with} $\Psi_C$ if
\begin{equation}\label{commutes}
(\id\tens\sigma_{B,C})\comp(\Psi_B\vt\Psi_C)=\Psi_C\vt\Psi_B,
\end{equation}
where $\sigma_{B,C}\in\Mor(B\tens{C},C\tens{B})$ is the flip.

This is a straightforward generalization of the notion of commutation of
classical families of maps. Note that $\Psi_B$ commutes with $\Psi_C$ if and
only if $\Psi_C$ commutes with $\Psi_B$.

The most basic properties of this notion will be
analyzed in the next proposition. In its formulation we use the concept of a
trivial quantum family of maps defined in Definition \ref{DefQFM}
\eqref{DefQFM4}.

\begin{prop}\label{comu}
Let $B$ be a $\cst$-algebra and let $\Psi_B\in\Mor(M,M\tens{B})$ be a quantum
family of maps $\qs(M)\to\qs(M)$ labeled by $\qs(B)$. Then
\begin{enumerate}
\item\label{comu1} if $\Psi_B$ is trivial then any quantum family of maps
$\qs(M)\to\qs(M)$ commutes $\Psi_B$.
\item\label{comu2} If $\Psi_B$ commutes with $\Phi\in\Mor(M,M\tens{A})$ then
$\Psi_B$ is trivial.
\end{enumerate}
\end{prop}

\begin{proof}
Ad \eqref{comu1}. Let $C$ be a $\cst$-algebra and let
$\Psi_C\in\Mor(M,M\tens{C})$ be a quantum family of maps $\qs(M)\to\qs(M)$.
Since $\Psi_B$ is assumed to be trivial, we have
\[
\begin{split}
(\id\tens\sigma_{C,B})(\Psi_C\vt\Psi_B)(m)
&=(\id\tens\sigma_{B,C})(\Psi_C\tens\id)(m\tens{I})\\
&=(\id\tens\sigma_{B,C})\bigl(\Psi_C(m)\tens{I}\bigr)\\
&=(\Psi_B\tens\id)\Psi_C(m)=(\Psi_C\vt\Psi_B)(m)
\end{split}
\]
for any $m\in{M}$.

Ad \eqref{comu2}. Let $\Lambda\in\Mor(A,M)$ be the unique morphism such that
$(\id\tens\Lambda)\Phi(m)=I\tens{m}$ for all $m\in{M}$ (see proofs of
Propositions \ref{ergodic} and \ref{erg-omega}). Let us apply
$(\id\tens\Lambda\tens\id)$ to both sides of
\[
(\id\tens\sigma_{B,A})(\Psi_B\vt\Phi)(m)=(\Phi\vt\Psi_B)(m)
\]
we obtain $I\tens{m}\tens{I}=I\tens\Psi_B(m)$. Since $m$ is arbitrary, $\Psi_B$
must be trivial.
\end{proof}

If two quantum families of maps $\qs(M)\to\qs(M)$ commute with a third one,
then so does their composition.

\begin{prop}\label{comcom}
Let $B$, $B'$ and $C$ be $\cst$-algebras and let
\[
\begin{split}
\Psi_B&\in\Mor(M,M\tens{B}),\\
\Psi_{B'}&\in\Mor(M,M\tens{B'}),\\
\Psi_C&\in\Mor(M,M\tens{C})
\end{split}
\]
be quantum families of maps $\qs(M)\to\qs(M)$ labeled by $\qs(B)$, $\qs(B')$
and $\qs(C)$ respectively. Assume that $\Psi_B$ commutes with $\Psi_C$ and that
$\Psi_{B'}$ commutes with $\Psi_C$. Then $\Psi_{B}\vt\Psi_{B'}$ commutes with
$\Psi_C$.
\end{prop}

\begin{proof}
In the following computation we shall use the fact that
\[
\sigma_{B\tens{B'},C}=(\sigma_{B,C}\tens\id)\comp(\id\tens\sigma_{B',C}).
\]
Now
\[
\begin{split}
(\id\tens\sigma_{B\tens{B'},C})\circ
&\bigl((\Psi_B\vt\Psi_{B'})\vt\Psi_C\bigr)
=(\id\tens\sigma_{B\tens{B'},C})
\comp(\Psi_B\vt\Psi_{B'}\vt\Psi_C)\\
&=(\id\tens\sigma_{B\tens{B'},C})\comp
\bigl((\Psi_B\tens\id\tens\id)\comp(\Psi_{B'}\tens\id)\comp\Psi_C\bigr)\\
&=(\id\tens\sigma_{B,C}\tens\id)\comp(\id\tens\id\tens\sigma_{B',C})\comp
(\Psi_B\tens\id\tens\id)\comp(\Psi_{B'}\tens\id)\comp\Psi_C\\
&=(\id\tens\sigma_{B,C}\tens\id)\comp(\Psi_B\tens\id\tens\id)\comp
(\id\tens\sigma_{B',C})\comp(\Psi_{B'}\tens\id)\comp\Psi_C\\
&=(\id\tens\sigma_{B,C}\tens\id)\comp(\Psi_B\tens\id\tens\id)\comp
(\id\tens\sigma_{B',C})\comp(\Psi_{B'}\vt\Psi_C)\\
&=(\id\tens\sigma_{B,C}\tens\id)\comp(\Psi_B\tens\id\tens\id)\comp
(\Psi_C\vt\Psi_{B'})\\
&=(\id\tens\sigma_{B,C}\tens\id)\comp(\Psi_B\tens\id\tens\id)\comp
(\Psi_C\tens\id)\comp\Psi_{B'}\\
&=\Bigl(\bigl[(\id\tens\sigma_{B,C})\comp(\Psi_B\tens\id)\comp\Psi_C\bigr]
\tens\id\Bigr)\comp\Psi_{B'}\\
&=\Bigl(\bigl[(\id\tens\sigma_{B,C})\comp(\Psi_B\vt\Psi_C)\bigr]
\tens\id\Bigr)\comp\Psi_{B'}\\
&=\bigl([\Psi_C\vt\Psi_B]\tens\id\bigr)\comp\Psi_{B'}\\
&=(\Psi_C\vt\Psi_B)\vt\Psi_{B'}=\Psi_C\vt(\Psi_B\vt\Psi_{B'})
\end{split}
\]
means that $\Psi_B\vt\Psi_{B'}$ commutes with $\Psi_C$.
\end{proof}

\begin{rem}
The notion of commuting quantum families of maps can be generalized further by
introducing braiding. The point is that instead of the flip $\sigma_{B,C}$ we
can put in \eqref{commutes} a braiding automorphism, say, $\Sigma_{B,C}$ (see
e.g.~\cite[Section 9.2]{majid}). Propositions \ref{comu} and \ref{comcom}
remain valid, but commutation fails to be a symmetric relation unless
$\Sigma_{B,C}^{-1}=\Sigma_{C,B}$
\end{rem}

Equipped with the conclusion of Proposition \ref{comcom} we can now easily
proceed to construct the quantum family of all maps $\qs(M)\to\qs(M)$ commuting
with a given quantum family. Then we find that it has a structure of a quantum
semigroup. The strategy is very similar to the one used in construction of
$\QMap^\omega\bigl(\qs(M)\bigr)$ given in Section \ref{invS}.

\begin{thm}\label{QcomPsi}
Let $B$ be a $\cst$-algebra and let $\Psi_B\in\Mor(M,M\tens{B})$ be a quantum
family of maps $\qs(M)\to\qs(M)$ labeled by $\qs(B)$. Let $K$ be the ideal of
$A$ generated by
\[
\Bigl\{(\omega\tens\id\tens\eta)(\Phi\vt\Psi_B)(m)
-(\omega\tens\eta\tens\id)(\Psi_B\vt\Phi)(m)
\st{m}\in{M},\:\omega\in{M^*},\:\eta\in{B^*}\Bigr\}.
\]
Let $\bar{A}$ be the quotient $A/K$, let $\rho:A\to\bar{A}$ be the canonical
epimorphism and let $\bar{\Phi}=(\id\tens\rho)\comp\Phi$. Then
\begin{enumerate}
\item the quantum family $\bar{\Phi}\in\Mor\bigl(M,M\tens\bar{A}\bigr)$ of maps
$\qs(M)\to\qs(M)$ labeled by $\qs\bigl(\bar{A}\bigr)$ commutes with the quantum
family $\Psi_B$.
\item\label{QcomPsi2} For any $\cst$-algebra $C$ and any quantum family
$\Psi_C\in\Mor(M,M\tens{C})$ of maps $\qs(M)\to\qs(M)$ which commutes with
$\Psi_B$ there exists a unique $\Lambda\in\Mor\bigl(\bar{A},C\bigr)$ such that
$\Psi_C=(\id\tens\Lambda)\comp\bar{\Phi}$.
\item There exists a unique
$\bar{\Delta}\in\Mor\bigl(\bar{A},\bar{A}\tens\bar{A}\bigr)$
such that
\[
\bigl(\bar{\Phi}\tens\id\bigr)\comp\bar{\Phi}=\bigl(\id\tens\bar{\Delta}
\bigr)\comp\bar{\Phi}.
\]
\item The morphism $\bar{\Delta}$ is coassociative and
$\bigl(\bar{A},\bar{\Delta}\bigr)$ is a compact quantum semigroup with unit.
$\bar{\Phi}$ is an action of $\bigl(\bar{A},\bar{\Delta}\bigr)$ on $\qs(M)$.
\item\label{QcomPsi4} $\rho$ is a quantum semigroup morphism.
\item\label{QcomPsi5} For any quantum semigroup $(C,\Delta_C)$ acting on
$\qs(M)$ with
action $\Phi_C\in\Mor(M,M\tens{C})$ commuting with $\Psi_B$, the unique
morphism $\Lambda\in\Mor\bigl(\bar{A},C\bigr)$ such that
$\Phi_C=(\id\tens\Lambda)\comp\bar{\Phi}$ is a quantum semigroup morphism.
\end{enumerate}
\end{thm}

In the proof we will only indicate the main steps. The details are practically
identical to those in the proofs of Theorems \ref{qstr} and \ref{Qw} as well as
Proposition \ref{qsm}.

\begin{proof}[Proof of Theorem \ref{QcomPsi}]
Clearly $\bar{\Phi}\in\Mor\bigl(M,M\tens\bar{A}\bigr)$ is a quantum family of
maps commuting with $\Psi_B$ and it is universal in the sense of Statement
\eqref{QcomPsi2}. Therefore, by Proposition \ref{comcom}, there exists a
unique morphism $\bar{\Delta}\in\Mor\bigl(\bar{A},\bar{A}\tens\bar{A}\bigr)$
such that
\[
\bigl(\id\tens\bar{\Delta}\bigr)\comp\bar{\Phi}=\bar{\Phi}\vt\bar{\Phi}
=\bigl(\bar{\Phi}\tens\id\bigr)\comp\bar{\Phi}.
\]

The coassociativity of $\bar{\Delta}$ follows exactly as in the proof of
Theorem \ref{qstr}\eqref{qstr2}. We use the fact that
\[
\bigl\{(\eta\tens\id)\bar\Phi(m)\st{m}\in{M},\:\eta\in{M^*}\bigr\}
\]
generates $\bar{A}$ ($\rho$ is an epimorphism). Thus
$\bigl(\bar{A},\bar{\Delta}\bigr)$ is a compact quantum semigroup with unit and
$\bar{\Phi}$ is its action on $\qs(M)$ which, by construction, commutes with
the quantum family of maps $\Psi_B$.

Statement \eqref{QcomPsi4} follows from Proposition \ref{qsm} and
\eqref{QcomPsi5} can be proved along the same lines as that proposition.
\end{proof}

The quantum semigroup $\bigl(\bar{A},\bar{\Delta}\bigr)$ constructed in Theorem
\ref{QcomPsi} will be denoted by the symbol $\QMap_{\Psi_B}\bigl(\qs(M)\bigr)$.

\section{Remarks on quantum groups inside quantum semigroups}\label{gr}

Let $M$ be, as in previous sections, a finite dimensional $\cst$-algebra. We
have constructed the compact quantum semigroup $\QMap\bigl(\qs(M)\bigr)$ and
have shown how to find quantum subsemigroups $\QMap^\omega\bigl(\qs(M)\bigr)$
for a state $\omega$ and $\QMap_{\Psi_B}\bigl(\qs(M)\bigr)$ for a given quantum
family $\Psi_B\in\Mor(M,M\tens{B})$ of maps $\qs(M)\to\qs(M)$ labeled by
$\qs(B)$. An interesting question arises: when are the quantum semigroups
$\QMap^\omega\bigl(\qs(M)\bigr)$ or $\QMap_{\Psi_B}\bigl(\qs(M)\bigr)$ compact
quantum groups?

Let us note that the quantum groups of automorphisms of of finite spaces
considered by S.~Wang in \cite{Qsym} are exactly the semigroups
$\QMap^\omega\bigl(\qs(M)\bigr)$ for $\omega$ the canonical trace (see
paragraph following the proof of Theorem \ref{grupa}).
Considering these examples one is quickly lead to a conjecture that
in general, $\QMap^\omega\bigl(\qs(M)\bigr)$ is a compact quantum group if
$\omega$ is faithful. Similarly one expects a compact quantum group structure
on $\QMap_{\Psi_B}\bigl(\qs(M)\bigr)$ for an ergodic family $\Psi_B$.

We will give here a short argument that a compact quantum semigroup acting on
$\qs(M)$, preserving a faithful state $\omega$ and satisfying a condition
corresponding to \eqref{exSLW2} form Theorem \ref{exSLW} (cf.~Remark 
\ref{ost}\eqref{ost3}), has cancellation from
the left (see e.g.~\cite[Remark 3]{cqg}). Moreover, one natural condition 
(cf.~Remark \ref{ost}\eqref{ost3})
immediateley guarantees that it also has cancellation from the right and is thus
a compact quantum group. We will first concentrate on the case when $\omega$ is a trace. This is the case encountered in the construction of S.~Wang's quantum automorphism groups of finite spaces (\cite{Qsym}) which we will present within the framework of quantum families of maps. Then we will present the result for a general faithul state.

Before we procede let us recall that if $(B,\Delta_B)$ is a compact quantum
semigroup then an \emph{$n$-dimensional representation} of $(B,\Delta_B)$ is an
$n\times{n}$-matrix $V=(v_{k,l})_{k,l=1,\ldots,n}$ of elements of $B$ such that
\begin{equation}\label{delVW}
\Delta_B(v_{k,l})=\sum_{r=1}^nv_{k,r}\tens{v_{r,l}}.
\end{equation}
We can also consider $V$ as an element in $M_n\tens{B}$ (where $M_n$ is the
$n\times{n}$-matrix algebra).

If $V\in{M_n}\tens{B}$ and $W\in{M_k}\tens{B}$ are representations of
$(B,\Delta_B)$ of dimensions $n$ and $k$ respectively, then the \emph{tensor
product} $V\tp{W}$ of $V$ and $W$ is the $nk$-dimensional representation
\[
V\tp{W}=V_{13}W_{23}\in{M_{nk}}\tens{B}.
\]
The matrix elements of $V\tp{W}$ are all possible products of a matrix element
of $V$ and a matrix element of $W$. Note also  that if $V$ and $W$ satisfy
$V^*V=I$ and $W^*W=I$ then $(V\tp{W})^*(V\tp{W})=I$. Similarly $VV^*=I$ and
$WW^*=I$ imply $(V\tp{W})(V\tp{W})^*=I$.

Our argument will on purpose be made in terms of martices over
$\cst$-algebras instead of adjointable maps of Hilbert $\cst$-modules. We hope
that this will make it more transparent to a reader not familiar with these
concepts.

\begin{lem}\label{lemV}
Let $(B,\Delta_B)$ be a quantum semigroup. Let
$V=(v_{k,l})_{k,l=1,\ldots,n}\in{M_n}\tens{B}$
be an $n$-dimensional representation of $(B,\Delta_B)$ such that $V^*V=I$. Then
for any $c\in{B}$ and any $k,l\in\{1,\ldots,n\}$ the element $c\tens{v_{k,l}}$
belongs to the linear span of the set
$\bigl\{(a\tens{I})\Delta_B(b)\st{a,b\in{B}}\bigr\}$.
\end{lem}

\begin{proof}
Equation \eqref{delVW} can be rewritten in matrix form as
\[
\begin{pmatrix}
\Delta_B(v_{1,1})&\cdots&\Delta_B(v_{1,n})\\
\vdots&&\vdots\\
\Delta_B(v_{n,1})&\cdots&\Delta_B(v_{n,n})
\end{pmatrix}
=\begin{pmatrix}
v_{1,1}\tens{I}&\cdots&v_{1,n}\tens{I}\\
\vdots&&\vdots\\
v_{n,1}\tens{I}&\cdots&v_{n,n}\tens{I}\\
\end{pmatrix}
\begin{pmatrix}
I\tens{v_{1,1}}&\cdots&I\tens{v_{1,n}}\\
\vdots&&\vdots\\
I\tens{v_{n,1}}&\cdots&I\tens{v_{n,n}}\\
\end{pmatrix}
\]
Multiplying this equality from the left by
\[
\begin{pmatrix}
{v_{1,1}^*}\tens{I}&\cdots&{v_{n,1}^*}\tens{I}\\
\vdots&&\vdots\\
{v_{1,n}^*}\tens{I}&\cdots&{v_{n,n}^*}\tens{I}\\
\end{pmatrix}
\]
we get all elements of the form $I\tens{v_{k,l}}$ as linear combinations of
elements of $(B\tens{I})\Delta_B(B)$. Multiplying by $c\tens{I}$ yields the
result.
\end{proof}

In the same way we can prove the following lemma:

\begin{lem}\label{lemVst}
Let $(B,\Delta_B)$ be a quantum semigroup. Let
$V=(v_{k,l})_{k,l=1,\ldots,n}\in{M_n}\tens{B}$
be an $n$-dimensional representation of $(B,\Delta_B)$ such that $VV^*=I$. Then
for any $c\in{B}$ and any $k,l\in\{1,\ldots,n\}$ the element $v_{k,l}\tens{c}$
belongs to the linear span of the set
$\bigl\{\Delta_B(a)(I\tens{b})\st{a,b\in{B}}\bigr\}$.
\end{lem}

\begin{thm}\label{grupa}
Assume that $\omega$ is a faithful trace and let $(B,\Delta_B)$ be a compact
quantum semigroup acting on $\qs(M)$ with action $\Phi_B$ and that
\begin{equation}\label{genB}
\Bigl\{(\eta\tens\id)\Phi_B(m)\st{m}\in{M},\:\eta\in{M^*}\Bigr\}
\end{equation}
generates $B$ as a $\cst$-algebra. Assume further that $\omega$ is invariant
for $\Phi_B$. Then $(B,\Delta_B)$ has cancellation from the left, i.e.
$\bigl\{(a\tens{I})\Delta_B(b)\st{a,b}\in{B}\bigr\}$
is linearly dense in $B\tens{B}$.

If, moreover, the set
\begin{equation}\label{comabyc}
\Bigl\{\Phi_B(m)(I\tens{a})\st{m}\in{M},\:a\in{B}\Bigr\}
\end{equation}
is linearly dense in $M\tens{B}$ then $(B,\Delta_B)$ has cancellation from the
right, i.e. the linear span
of $\bigl\{\Delta_B(a)(I\tens{b})\st{a,b}\in{B}\bigr\}$
is dense in $B\tens{B}$.
\end{thm}

\begin{proof}
Let $(m_l)_{l=1,\ldots,n}$ be a basis of $M$ which is orthonormal for the scalar
product given by $\omega$. Let $\widetilde{a}=(a_{k,l})_{k,l=1,\ldots,n}$ be
the matrix of elements of $B$ such that
\begin{equation}\label{defaij}
\Phi_B(m_l)=\sum_{k=1}^nm_k\tens{a_{k,l}}.
\end{equation}
Then $\widetilde{a}$ is an $n$-dimensional representation of $(B,\Delta_B)$.

The invariance of $\omega$ for $\Phi_B$ implies that $\widetilde{a}$ is an
isometry:
\begin{equation}\label{compu}
\begin{split}
\delta_{l,t}I=\omega(m_l^*m_t)I
&=(\omega\tens\id)\Phi_B(m_l^*m_t)
=(\omega\tens\id)\bigl(\Phi_B(m_l)^*\Phi_B(m_t)\bigr)\\
&=(\omega\tens\id)\left(\Bigl(\sum_km_k^*\tens{a_{k,l}^*}\Bigr)
\Bigl(\sum_sm_s\tens{a_{s,t}}\Bigr)\right)\\
&=(\omega\tens\id)\Bigl(\sum_{k,s}m_k^*m_s\tens{a_{k,l}^*}a_{s,t}\Bigr)\\
&=\sum_{k,s}\omega(m_k^*m_s){a_{k,l}^*}a_{s,t}
=\sum_{k,s}\delta_{k,s}{a_{k,l}^*}a_{s,t}\\
&=\sum_k{a_{k,l}^*}a_{k,t}=\bigl[\widetilde{a}^*\widetilde{a}\bigr]_{l,t}
\end{split}
\end{equation}

The fact that $\omega$ is a trace means that $(m_l^*)_{l=1,\ldots,n}$ is also
an orthonormal basis of $M$ for the scalar product given by $\omega$. Therefore
the matrix $\widetilde{b}$ with $a_{k,l}^*$ as the $(k,l)$-entry, is also an
$n$-dimensional representation of $(B,\Delta_B)$ and an isometry.

Let $V$ be a representation of $(B,\Delta_B)$ constructed as a finite
tensor product of representations $\widetilde{a}$ and $\widetilde{b}$ (in any
order). Then $V$ satisfies the conditions in Lemma \ref{lemV} and so, for any
matrix element $x$ of $V$ and any $c\in{B}$, the element $c\tens{x}$ belongs to
$(B\tens{I})\Delta_B(B)$.

Clearly any monomial in matrix elements of $\widetilde{a}$ and $\widetilde{b}$
is a matrix element of $V$ of considered form. Our assumption that $B$ is
generated by \eqref{genB} can be rephrased as saying that the span of all
monomials in matrix elements of $\widetilde{a}$ and $\widetilde{b}$ is dense in
$B$. This ends the proof that $(B,\Delta_B)$ has cancellation from the left.

Now let us note that the closure of the span of \eqref{comabyc} is the image of
$\widetilde{a}$ considered as a map of $M\tens{B}\cong{B^n}$ into itself. Thus from the
linear density of \eqref{comabyc} in $M\tens{B}$ it follows that
$\widetilde{a}$ is unitary and we might as well apply our reasoning to
$\widetilde{a}^*$. Using Lemma \ref{lemVst} we arrive at the conclusion that
$(B,\Delta_B)$ as cancellation from the right.
\end{proof}

\begin{rem}
Condition that \eqref{genB} generated $B$ of Theorem \ref{grupa} is there to ensure that the action of $(B,\Delta_B)$ on $\qs(M)$ is in some sense faithful. One can imagine a large quantum semigroup acting via a quotient group and \eqref{genB} excludes such a situation. The condition of density of \ref{comabyc} has been used before in connection with quantum group actions. It has been introduced by P.~Podle\'s in his thesis \cite[Definicja 2.2]{podPHD}.
\end{rem}

\newcommand{\AutX}{\mathrm{Aut}(X_n)}

Let us now now see how the quantum permutation groups of Wang (\cite{Qsym}) 
fit into the framework described above and in Section \ref{invS}. In 
\cite{Qsym} S.~Wang described quantum automorphism groups of finite spaces. 
For each natural $n$ Wang considered the finite space $X_n=\{1,\ldots,n\}=\qs(\CC^n)$ 
and defined the quantum group $\AutX=(A,\Delta)$ where $A$ is the 
$\cst$-algebra generated by elements 
\begin{equation}\label{gener}
\bigl\{a_{ij}\st{i,j=1,\ldots,n}\bigr\}
\end{equation}
with relations
\begin{subequations}\label{WangR}
\begin{align}
(a_{ij})^2&=a_{ij},&i,j&=1,\ldots,n,\label{Wang1}\\
(a_{ij})^*&=a_{ij},&i,j&=1,\ldots,n,\label{Wang2}\\
\sum_{j=1}^na_{ij}&=I,&i&=1,\ldots,n,\label{Wang3}\\
\sum_{i=1}^na_{ij}&=I,&j&=1,\ldots,n,\label{Wang4}
\end{align}
\end{subequations}
It turns out that 
\begin{equation}\label{row}
\AutX=\QMap^{\psi}(X_n), 
\end{equation}
where
$\psi$ is the uniformly distributed probability measure on $X_n$. Indeed,
if we let $\widetilde{A}$ be the $\cst$-algebra generated by \eqref{gener} 
with relations \eqref{Wang1}, \eqref{Wang2} and \eqref{Wang3} then
$\qs\bigl(\widetilde{A}\bigr)$ is the quantum space of all maps 
$X_n\to{X_n}$. The comultiplication 
\[
\widetilde{\Delta}:a_{ij}\longmapsto\sum_{k=1}^na_{ik}\tens{a_{kj}}
\]
coincides with the one described in Theorem \ref{qstr}\eqref{qstr1}. In other words
\[
\bigl(\widetilde{A},\widetilde{\Delta}\bigr)=\QMap(X_n).
\]
Now the relation \eqref{Wang4} defines the ideal $J$ related to $\psi$ 
as in Theorem \ref{Qw} and the comultiplication $\Delta$ on $A$ is given by 
the same formula as $\widetilde{\Delta}$. All this shows that \eqref{row} holds.

The action $\Phi$ of $(A,\Delta)$ on $X_n$ is given by
\[
\Phi(e_j)=\sum_{i=1}^ne_i\tens{a_{ij}},
\]
where $\{e_1,\ldots,e_n\}$ is the standard basis of $\CC^n=\C(X_n)$. Now we will
show that $\QMap^\psi(X_n)$ is in fact a compact quantum \emph{group} and not
merely a quantum semigroup. First let us note that the density condition \eqref{genB} is 
satisfied by Remark \ref{ost}\eqref{ost3} (with $M=\CC^n$). 
Moreover linear density of \eqref{comabyc}
in $\CC^n\tens{A}$ also holds. This follows from the fact that we have\footnote{
Any family $\{p_1,\ldots,p_n\}$ of projections such that $\sum\limits_{k=1}^np_k=I$ must satisfy
$p_kp_l=\delta_{kl}p_k$ for $k,l=1,\ldots,n$.} 
\[
a_{ij}a_{kj}=\delta_{ik}a_{kj}
\]
and consequently 
\[
\begin{split}
\sum_{j=1}^n\Phi(e_j)(I\tens{a_{kj}})
&=\sum_{j=1}^n\sum_{i=1}^ne_i\tens{a_{ij}a_{kj}}\\
&=\sum_{j=1}^n\sum_{i=1}^ne_i\tens\delta_{ik}a_{kj}\\
&=\sum_{i=1}^ne_i\tens\biggl(\delta_{ik}\sum_{j=1}^na_{kj}\biggr)\\
&=e_k\tens\biggl(\sum_{j=1}^na_{kj}\biggr)=e_k\tens{I}.
\end{split}
\]
Therefore $m\tens{a}\in$ belongs to the span of \eqref{comabyc} for any $m\in{\CC^n}$ and 
any $a\in{A}$. 

All this shows that the assumptions of Theorem \ref{grupa} are satisfied and
$(A,\Delta)=\QMap^\psi(X_n)$ is a compact quantum group. Let us emphasize that 
this has already been established in \cite{Qsym}. The point of the above argument
is that the existence of a quantum group structure on this quantum space is a 
consequence of its more fundamental properties, namely that $\qs(A)$ is
the universal quantum family of maps $X_n\to{X_n}$ preserving the measure $\psi$.

The statement of Theorem \ref{grupa} gives explicit relationship between cancellation laws and density conditions. However the conclusion that $(B,\Delta_B)$ is a compact quantum group does not require that $\omega$ be a trace. In facr we have

\begin{thm}
Let $(B,\Delta_B)$ be a quantum semigroup and let $\Phi_B\in\Mor(M,M\tens{B})$ be an action of $(B,\Delta_B)$ on a finite dimensional $\cst$-algebra $M$. Let $\omega$ be an ivvarian faithful state on $M$. Then if
\[
\Bigl\{(\eta\tens\id)\Phi_B(m)\st{m}\in{M},\:\eta\in{M^*}\Bigr\}
\]
generates $B$ as a $\cst$-algebra and the set
\begin{equation}\label{podle}
\Bigl\{\Phi_B(m)(I\tens{a})\st{m}\in{M},\:a\in{B}\Bigr\}
\end{equation}
is linearly dense in $M\tens{B}$ then $(B,\Delta_B)$ is a compact quantum group.
\end{thm}

The proof is almost the same as that of Theorem \ref{grupa} with the difference that the basis $(m_l^*)_{l=1,\ldots,n}$ is not orthonormal. However the matrix $\widetilde{a}$ with matrix element $(a_{i,j})_{i,j=1\ldots,n}$ defined by \eqref{defaij} still is an isometry. The density condition \eqref{podle} guarantees that $\widetilde{a}$ is unitary. Let $\Bar{\widetilde{a}}$ be the matrix with elements $(a_{i,j})_{i,j=1\ldots,n}$. Then $\Bar{\widetilde{a}}$ is not isometric. Still, by elementary linear algebra, there exists an invertible map $\sigma:M\to{M}$ such that
\[
\omega(yx)=\omega\bigl(y\sigma(x)\bigr)
\]
for all $x,y\in{M}$. Let $S=(s_{i,j})_{i,j=1,\ldots,n}$ be the matrix of $\sigma$ in the basis $(m_l)_{l=1,\ldots,n}$. 
\[
\omega(m_im_j^*)=\omega\bigl(m_j^*\sigma(m_i)\bigr)=\sum_{p=1}^ns_{i,p}\omega(m_j^*m_p)=s_{i,j}.
\]
Therefore the computation \eqref{compu} shows that
\[
s_{i,j}I=\omega(m_im_j^*)I=\sum_{p,q=1}^n\omega(m_pm_q^*)a_{p,i}a_{q,j}^*
=\sum_{p,q=1}^na_{p,i}s_{p,q}a_{q,j}^*
=\Bigl[{\Bar{\widetilde{a}}}^*(S\tens{I})\Bar{\widetilde{a}}\Bigr]_{i,j}
\]
or in other words
\begin{equation}\label{sij}
S\tens{I}={\Bar{\widetilde{a}}}^*(S\tens{I})\Bar{\widetilde{a}}.
\end{equation}
Multiplying \eqref{sij} from the left by $S^{-1}\tens{I}$ we see that $\Bar{\widetilde{a}}$ is left invertible. However, as a map $B^n\to{B^n}$ the matrix $\Bar{\widetilde{a}}$ has dense range because of the density condition \eqref{podle}. It follows that $\Bar{\widetilde{a}}$ is invertible, and thus so is
\[
\widetilde{a}^\top=\Bar{\widetilde{a}}^*.
\]
Threfore $B$ is a $\cst$-algebra generated by elements of a unitary matrix $\widetilde{a}=(a_{i,j})_{i,j=1,\ldots,n}$ with a morphism $\Delta_B\in\Mor(B,B\tens{B})$ satisfying
\[
\Delta_B(a_{i,j})=\sum_{p=1}^na_{i,p}\tens{a_{p,j}}
\]
(this follows from \eqref{defaij} and the fact that $\Phi_B$ is an action) and such that $\widetilde{a}^\top$ is invertible. By the results of \cite{remqg} $(B,\Delta_B)$ is a compact quantm group.

\end{document}